\documentclass[12pt, a4paper, reqno, oneside]{amsart}
\usepackage [utf8]{inputenc}
\usepackage{mathtools}
\usepackage{amsthm}
\usepackage{amssymb}

\usepackage{cite, verbatim}

\newtheorem{Theorem}{Theorem}
\newtheorem{cor}{Corollary}
\newtheorem{question}{Problem}

\begin{document}

\title{On $p$-index extremal groups}

\author{Andrey V. Vasil'ev and Ilya B. Gorshkov}

\thanks{A.\,V.\,Vasil'ev was supported by National Natural Science Foundation of China (No.~12171126), I. B. Gorshkov was supported by the Russian Science Foundation, project  18-71-10007.}

\begin{abstract} The question on connection between the structure of a finite group $G$ and the properties of the indices of elements of $G$ has been a popular research topic for many years. The {\em$p$-index $|x^G|_p$} of an element $x$ of a group $G$ is the $p$-part of its index~$|x^G|=|G:C_G(x)|$. The presented short note describes some new results and open problems in this direction, united by the concept of the $p$-index of a group element.
\smallskip

\noindent\textsc{Key words:} finite group, conjugacy class, index, $p$-index of a group element.

\noindent\textsc{MSC:} 20E45, 20D60.
 \end{abstract}
 
\maketitle

The size of the conjugacy class of an element of a finite group, which is known to be equal to the index of its centralizer, is one of main numerical characteristics of a finite group. The question on connection between the structure of a finite group and the properties of the size set of its conjugacy classes has been a popular research topic for many years. Thus, at the very beginning of the 20th century, W.\,Burnside  in~\cite{Burnside} showed  that a group of composite order, having a conjugacy class whose size is a power of a prime number, cannot be simple. That assertion was the key to his famous theorem on the solvability of finite $\{p,q\}$-groups. The presented short note describes some new results and open problems in this direction, united by the concept of the {\em$p$-index} of a group element.

First, we introduce  necessary definitions and notation. Let $G$ be a finite group and $x\in G$. As usual, $Z(G)$ denotes the center of~$G$, $C_G(x)$ denotes the centralizer of $x$ in~$G$, and $x^G$ denotes the conjugacy class in~$G$ containing the element~$x$, in other words, $x^G$ is the orbit of $x$ with respect to the action by conjugacy of $G$ on itself. The size of the given orbit (class) $|x^G|=|G:C_G(x)|$ is called the {\em index} of the element $x$ in $G$. The set of all indices of the elements of $G$ is denoted by $N(G)$ (the notation $\operatorname{cl}(G)$ is also very common, see, e.\,g., \cite{CasTom}).

If $p$ is a prime ($\pi$ is some set of primes), then for a positive integer $n$, we refer to $n_p$ as the $p$-part of $n$, i.\,e. the largest power of $p$ dividing $n$ (respectively, to $n_\pi$ as the $\pi$-part of~$n$). The {\em$p$-index $|x^G|_p$} of an element $x$ of a group $G$ is the $p$-part of~$|x^G|$. For a group $G$, we set $N(G)_p=\{|x^G|_p : x\in G\}$ and call the {\em$p$-index $|G||_p$} of $G$ the largest element in the set $N(G)_p$, i.\,e. $|G||_p=\max\{|x^G|_p : x\in G\}$.

Following \cite{CasTom}, we denote by $e_p(G)$ the number of nontrivial elements in~$N(G)_p$. We are interested in groups $G$ with $e_p(G)=1$; we call them \emph{extremal with respect to} $p$-\emph{index} or, simply, {\em $p$-index extremal} and write $G\in R(p)$ (in the recent paper by the second author \cite{GorRp} these groups were called {\em$R(p)$-groups}).

Although the class $R(p)$ has been singled out recently, results on groups in this class (with various additional restrictions) are well known. Already in 1953, N.\,Ito proved that in a group $G$ with $N(G)=\{1,n\}$ the number~$n$ is a power of some prime~$p$, and the group itself is nilpotent~\cite{Ito53}. More precisely, $G$ is the direct product of a $p$-group $A$ and an abelian group $B$ of order coprime to~$p$. It was later established in \cite{Ish} that the nilpotency class of the group $A$ does not exceed~$3$. In \cite{Camin72}, A.\,Kamina showed that for any primes $p$ and $q$ the group $G$ with $N(G)=\{1, p^\alpha, q^\beta, p^\alpha q^\beta\}$ is the direct product of a $p$-group $A$ and $q$-group $B$ and, in particular, is nilpotent, while Beltran and Felipe \cite{BelFel, BelFel2} generalized this result for the case of two arbitrary coprime positive integers $m$ and~$n$. Obviously, all the above groups will be extremal with respect to any prime $p$ (later we will return to groups with this property).

Recall that a group is called $p$-solvable if it contains a normal series whose factors are either $p$-groups or $p'$-groups. In \cite{CasTom}, Casolo and Tombari considered the $p$-index extremal $p$-solvable groups. It was shown that the $p$-length of every such group must be equal to~$1$. Moreover, if a Sylow $p$-subgroup $P$ of such group is nonabelian, then $G$ has a normal $p$-complement \cite[Theorem~2.3]{CasTom} and $Z(P)\leq Z(G) $ \cite[Lemma~2.4]{CasTom}. It turns out that the last two statements are true in a much more general situation. Note first that the class $R(p)$ splits naturally into two non-overlapping subclasses:

1) $G\in R(p)^*$, if there is a $p$-element in $G$ with a nontrivial $p$-index;

2) $G\in R(p)^{**}$, otherwise.

Recently, the second author obtained the following results.

\begin{Theorem}\label{th:rp*}{\em\cite[Theorem~1.3]{GorRp}}
If $G\in R(p)^*$, then $G$ has the normal $p$-complement.
\end{Theorem}

\begin{cor}\label{cor:rp*}{\em\cite[Corollary~1.4]{GorRp}}
If $G\in R(p)^{*}$ and $P$ is a Sylow $p$-subgroup of~$G$, then $Z(P)\leq Z(G)$.
\end{cor}

Certainly, if there is a $p$-element in $G$ with nontrivial $p$-index, then a Sylow $p$-subgroup $P$ is nonabelian. The converse is not true in general \cite[Theorem~B]{NavSolTiep}. Nevertheless, it is easy to deduce from \cite[Lemma~7]{GorGr} that in the case of a $p$-index extremal group, noncommutativity of a Sylow $p$-subgroup of $G$ yields $G\in R(p)^*$.

In \cite{Vas}, the first author noted that in a group $G$ with trivial center and $|G||_p=p$, a Sylow $p$-subgroup must be abelian (and even elementary abelian). Using Corollary~\ref{cor:rp*} in the case of $G\in R(p)^*$ and \cite[Theorem B]{NavTiep} in the case of $G\in R(p)^{**}$, it is easy to establish a similar assertion for an arbitrary $p$-index extremal group $G$ with trivial center.

\begin{cor}\label{cor:rpz=1}{\em\cite[Corollary~1.5]{GorRp}}
If $G\in R(p)$ and $Z(G)=1$, then a Sylow $p$-subgroup of $G$ is abelian.
\end{cor}

Although an $R(p)$-group $G$ with trivial center can have a noncyclic Sylow $p$-subgroup (for example, the group $A_5\in R(2)^{**}$ has this property), it seems that such groups can be fully classified.

\begin{question}\label{q:rp}
Describe the structure of the $p$-index extremal groups with trivial center and noncyclic Sylow $p$-subgroups.
\end{question}

Let us return to Theorem~\ref{th:rp*}. According to it, the group $G$ is a semidirect product of the normal $p'$-subgroup of $G$ and a $p$-subgroup~$P$. The structure of $P$ is generally clear (see \cite{Ish, Naik}), it would be interesting to have an answer to the following question.

\begin{question}\label{q:rp*}
Describe the structure of the normal $p$-complement of $G\in R(p)^*$.
\end{question}

It is easy to see that groups with the $R(p)^{**}$ property are quite common. For example, using the well-known properties of centralizers, one can show that for any finite nonabelian simple group $G$ there exists a prime $p$ such that a Sylow $p$-subgroup $S$ of $G$ is cyclic and $C_G(S)=C_G(x)$ for any nontrivial element $x\in S$, and hence $G\in R(p)^{**}$. It turns out that any group from the class $R(p)^{**}$ has the following property.

\begin{Theorem}\label{th:rp**}{\em\cite[Lemma~8]{GorGr}}
If $G\in R(p)^{**}$, then $G$ contains at most one nonabelian composition factor $S$ whose order is divisible by~$p$.
\end{Theorem}

We believe that more can be said in this situation.

\begin{question}\label{q:rp**}
Let $G\in R(p)^{**}$ be a nonsolvable group. Is it true that $G=P\times H$ is the direct product of a $p$-group $P$ and a group $H$ such that $O^{p'}(H/O_{p'}(H))$ is a nonabelian simple group?
\end{question}

The aforementioned results of Ito, Kamina and others dealt with groups that are extremal with respect to any $p$-index, we call them {\em index extremal} ({\em$\mathcal{D}$-groups} in the terminology from \cite{ CasTom}). These groups can also be described as groups with the following property:

\begin{equation}\label{eq:D}
N(G)\subseteq \prod_{p\in\pi(N(G))} N(G)_p,\mbox{ where }N(G)_p=\{1,p^{\alpha(p)}\}~\forall p\in\pi(N(G)).
\end{equation}

In the case when equality takes place in \eqref{eq:D} instead of inclusion, the group $G$ will be nilpotent. This fact was established in \cite[Theorem~1]{CasTom}. Moreover, in this paper, an explicit description was given for the index extremal groups satisfying the additional condition $p^{\alpha(p)}\in N(G)$ for any $p\in\pi(N(G))$ \cite[Theorem~2]{CasTom}. In particular, all such groups turn out to be solvable. In the general case, the latter is not true --- the nonabelian simple group $A_5$ is obviously index extremal. Therefore, the following problem appears.

\begin{question}\label{q:rpall}
Describe the structure of the finite index extremal groups.
\end{question}

Let $p$ and $q$ be distinct primes. A finite group is said to {\em have property} $\{p,q\}^*$ if $G\in R(p)\cap R(q)$ and $a_{\{p,q\}}>1$ for any $a\in N(G)\setminus\{1\}$. In \cite{Ito70}, Ito studied the simplest nontrivial subclass of $\{p,q\}^*$-groups, showing that if $N(G)=\{1,m,n\}$, where $m$ and $n$ are coprime numbers, then $G$ is solvable. In the general case, however, a $\{p,q\}^*$-group need not be solvable. Moreover, it can be shown (see, for example, \cite[Lemma~3.4(v)]{VasGr} and \cite[Lemma~3.2(v)]{GrZv} for classical groups) that, with some exceptions, for a simple group $G$ of Lie type, there are primes $p$ and $q$ such that $G$ is a $\{p,q\}^*$-group.

Recall that $|G||_p$ denotes the $p$-index of the group $G$. If $\pi$ is a set of primes, then we put $|G||_{\pi}=\prod_{p\in\pi}|G||_p$. It is obvious that $|G||_{\pi}$ divides $|G|_\pi$ for each~$\pi$. In \cite{Gorpq}, the following properties of $\{p,q\}^*$-groups with trivial center were obtained.

\begin{Theorem}\label{th:pq*}
Let $p$ and $q$ be distinct primes greater than $5$. If $G$ is a $\{p,q\}^*$-group with trivial center, then $|G|_{\{p,q\}}=|G||_{\{p, q\}}$.
\end{Theorem}

\begin{cor}\label{cor:pq*}
Let $p$ and $q$ be distinct primes greater than $5$. If $G$ is a $\{p,q\}^*$-group with trivial center, then $C_G(g)\cap C_G(h)=1$ for each nontrivial $p$-element $g$ and each nontrivial $q$-element $h$ of~$G$.
\end{cor}

The last corollary turned out to be the key one in the proof of an assertion about orthogonal groups of dimensions $8$ and $16$ \cite[Theorem~1]{GorThom}, which became the final step (see details in \cite{GorThom}) in proving the validity of the well-known conjecture stated by J.\,G.\,Thompson in a letter to W.\,Shi in 1987, and which in 1992 was puted by A.\,S.\,Kondratiev in the {\em Kourovka notebook}.

\medskip

\noindent\textbf{Thompson's Conjecture} \cite[Problem~12.38]{Kour}. {\it If $L$ is a finite nonabelian simple group, $G$ is a finite group with trivial center, and $N(G)=N(L)$, then $G$ and $L$ are isomorphic.}

\medskip

Certainly, the assertion of Thompson's conjecture cannot be extended to the case when $L$ is an arbitrary finite group with trivial center; in particular, there are a nonsolvable group $L$ and solvable group $G$, both with trivial center, satisfying $N(L)=N(G)$ \cite{Navv}. However, the following problem remains open.

\begin{question}\label{q:ThPower}{\em\cite[Problem~20.29]{Kour}}
Is it true that for any natural number $n$ and any nonabelian simple group $L$, if $G$ is a finite group with trivial center and $N(G)=N(P)$, where $P=L^n$ is the $n$th direct power of~$L$, then $G$ and $P$ are isomorphic?
\end{question}

In the case $n=1$, the positive answer to this question follows from the validity of Thompson's conjecture. If $n>1$, then at the moment there are only a few very particular results \cite{GorA5A5,PanA6A6,GorPan}.

\bigskip

Andrey V. Vasil'ev,

Sobolev Institute of Mathematics, Novosibirsk, RUSSIA;

School of Science, Hainan University, Haikou, Hainan, P. R. CHINA.

e-mail: vasand@math.nsc.ru

\medskip

Ilya B. Gorshkov,

Sobolev Institute of Mathematics, Novosibirsk, RUSSIA;

Siberian Federal University, Krasnoyarsk, RUSSIA.

e-mail: ilygor8@gmail.com

\end{document}